\documentclass[12pt]{amsart}
\usepackage{amssymb}

\begin{document}
\bibliographystyle{alpha}
\numberwithin{equation}{section}

\def\Label#1{\label{#1}}

\def\1#1{\ov{#1}}
\def\2#1{\widetilde{#1}}
\def\3#1{\mathcal{#1}}
\def\4#1{\widehat{#1}}

\def\s{s}
\def\k{\kappa}
\def\ov{\overline}
\def\span{\text{\rm span}}
\def\tr{\text{\rm tr}}
\def\GL{{\sf GL}}
\def\xo {{x_0}}
\def\Rk{\text{\rm Rk\,}}
\def\sg{\sigma}

\def \hn {holomorphically nondegenerate}
\def\hyp{hypersurface}
\def\prt#1{{\partial \over\partial #1}}
\def\det{{\text{\rm det}}}
\def\wob{{w\over B(z)}}
\def\co{\chi_1}
\def\po{p_0}
\def\fb {\bar f}
\def\gb {\bar g}
\def\Fb {\ov F}
\def\Gb {\ov G}
\def\Hb {\ov H}
\def\zb {\bar z}
\def\wb {\bar w}
\def \qb {\bar Q}
\def \t {\tau}
\def\z{\chi}
\def\w{\tau}
\def\Z{\zeta}

\def \T {\theta}
\def \Th {\Theta}
\def \L {\Lambda}
\def\b{\beta}
\def\a{\alpha}
\def\o{\omega}
\def\l{\lambda}

\def \im{\text{\rm Im }}
\def \re{\text{\rm Re }}
\def \Char{\text{\rm Char }}
\def \supp{\text{\rm supp }}
\def \codim{\text{\rm codim }}
\def \Ht{\text{\rm ht }}
\def \Dt{\text{\rm dt }}
\def \hO{\widehat{\mathcal O}}
\def \cl{\text{\rm cl }}
\def \bR{\mathbb R}
\def \bC{\mathbb C}
\def \bP{\mathbb P}
\def \C{\mathbb C}
\def \bL{\mathbb L}
\def \bZ{\mathbb Z}
\def \bN{\mathbb N}
\def \scrF{\mathcal F}
\def \scrK{\mathcal K}
\def \scrM{\mathcal M}
\def \cR{\mathcal R}
\def \scrJ{\mathcal J}
\def \scrA{\mathcal A}
\def \scrO{\mathcal O}
\def \scrV{\mathcal V}
\def \scrL{\mathcal L}
\def \scrE{\mathcal E}
\def \hol{\text{\rm hol}}
\def \aut{\text{\rm aut}}
\def \Aut{\text{\rm Aut}}
\def \J{\text{\rm Jac}}
\def\jet#1#2{J^{#1}_{#2}}
\def\gp#1{G^{#1}}
\def\gpo{\gp {2k_0}_0}
\def\emmp {\scrF(M,p;M',p')}
\def\rk{\text{\rm rk}}
\def\Orb{\text{\rm Orb\,}}
\def\Exp{\text{\rm Exp\,}}
\def\ess{\text{\rm Ess\,}}
\def\mult{\text{\rm mult\,}}
\def\Jac{\text{\rm Jac\,}}
\def\Span{\text{\rm span\,}}
\def\d{\partial}
\def\D{\3J}
\def\pr{{\rm pr}}
\def\dbl{[\![}
\def\dbr{]\!]}
\def\nl{|\!|}
\def\nr{|\!|}

\def \depth{\text{\rm depth\,}}
\def \D{\text{\rm Der}\,}
\def \Rk{\text{\rm Rk}\,}
\def \ima{\text{\rm im}\,}
\def \vfi{\varphi}

\title [Projection 
on Segre varieties]{Projection on Segre varieties and 
determination of holomorphic mappings between real submanifolds}
\author[M.S. Baouendi, P. Ebenfelt and L. P. Rothschild]{M.S. Baouendi, Peter Ebenfelt
and Linda P. Rothschild} \thanks{{\rm The first and third authors are supported in part by the NSF grant
DMS-0400880. The second
author is supported in part by the NSF grant DMS-0401215.
\newline}}

\dedicatory{Dedicated to Professor Sheng Gong on the occasion of 
his 75th birthday}

\address{ Department of Mathematics, University of
California at San Diego, La Jolla, CA 92093-0112, USA}
\email{sbaouendi@ucsd.edu, 
pebenfel@math.ucsd.edu,
lrothschild@ucsd.edu }


\abstract It is shown that a germ of a holomorphic mapping sending a real-analytic
generic submanifold of finite type into another is determined by its projection
on the Segre variety of the target manifold. A necessary and sufficient
condition is given for a germ of a mapping into the Segre variety of the target
manifold to be the projection of a holomorphic mapping sending the source
manifold into the target. An application to the biholomorphic equivalence
problem is also given. 
\endabstract

\newtheorem{Thm}{Theorem}[section]
\newtheorem{Def}[Thm]{Definition}
\newtheorem{Cor}[Thm]{Corollary}
\newtheorem{Pro}[Thm]{Proposition}
\newtheorem{Lem}[Thm]{Lemma}
\newtheorem{Rem}[Thm]{Remark}
\newtheorem{Ex}[Thm]{Example}

\maketitle
\section{Introduction and main results}  \Label{s:intro}

In this paper, we show that a germ at $0$ of a holomorphic mapping
$H\colon (\bC^N,0)\to (\bC^{N'},0)$ sending one real-analytic generic submanifold
$M\subset \bC^N$ of finite type at $0$ into another real-analytic generic submanifold
$M'\subset \bC^{N'}$ is determined by its projection onto $\Sigma'_{0}$, the
Segre variety of $M'$ at $0$ (Theorem \ref{T:Hdet1}). We also give a necessary
and sufficient condition for a germ at
$0$ of a holomorphic mapping
$F\colon (\bC^N,0)\to (\Sigma'_0,0)$  to be the projection of such a mapping $H$
(Theorem \ref{T:F}). As a corollary, we obtain a new criterion for two real-analytic
generic submanifolds $M,M'$ of finite type at $0\in \bC^N$ to be locally
biholomorphically equivalent at $0$ (Corollary \ref{C:F}). The main tools used in the
proofs are the iterated Segre mappings, as previously introduced by the authors in
\cite{BERalg}, \cite {BERJAG}, and also a new invariant description of normal
coordinates (Theorem
\ref{t:normal1}), which may be of independent interest. Segre variety
techniques in the context of mappings between real hypersurfaces were introduced
in \cite{webalg} and \cite{webrefl}.

 Let $M$ be a real-analytic generic submanifold of
codimension $d$ in $\bC^N$ with $0\in M$, given locally near
$0$   by
\begin{equation}\Label{e:defeq}
\rho_1(Z,\bar Z)=\ldots=\rho_d(Z,\bar Z)=0,\end{equation}
where $\rho(Z,\zeta):=(\rho_1(Z,\zeta),\ldots,\rho_d(Z,\zeta))$  is a $\bC^d$-valued
holomorphic function such that
$\partial_Z\rho_1\wedge\ldots\wedge\partial_Z\rho_d\neq 0$ near $0$ and
\begin{equation}\Label{e:reality}
\rho(Z,\zeta)=\overline{\rho(\bar \zeta,\bar Z)}.
\end{equation} 
The generic submanifold $M$ is said to be of {\it
finite type} at $p$ (in the sense of Kohn \cite{kohn} and
Bloom-Graham \cite{bgtype}) if the (complex) Lie algebra
$\frak g_M$ generated by all smooth $(1,0)$ and $(0,1)$
vector fields tangent to $M$ satisfies $\frak g_M(p)=\bC
T_{p}M$, where $\C T_{p} M$ is the complexified
tangent space to $M$ at $p$.

Let $U\subset \bC^N$ be a sufficiently small open neighborhood of $0$. 
For $p\in \bC^N$ sufficiently close to 0, we denote by $\Sigma_p$ the Segre variety of
$M$ at $p$ defined by 
\begin{equation}
\Sigma_p:=\{Z\in U\colon \rho(Z,\bar p)=0\}.
\end{equation}
We observe, for future reference, that $\Sigma_p$ is an $n$-dimensional complex
submanifold of
$U$, with
$n=N-d$, for all such $p$. Moreover, it follows from \eqref{e:reality} that $p\in
\Sigma_p$ if and only
$p$ is in $M$, and also that $q\in \Sigma_p$ if and only if $p\in \Sigma_q$.

\begin{Thm} \Label{T:Hdet1} Let $M$ be a real-analytic generic
submanifold of codimension $d$ in $\bC^N$ and of finite type at
$0\in M$. Then, for every $\lambda\in\bC$, 
$0<|\lambda|<1$,  there exist
$2d+1$ germs at
$0$ of holomorphic functions 
$g^\lambda_1,\ldots,g^\lambda_{d}$, $h^\lambda_1,\ldots, h^\lambda_{d+1}\colon
(\bC^N,0)\to \bC^N$, depending holomorphically on $\lambda$, such that
$g_j^{\lambda}(0)\to 0$, $h_j^\lambda(0)\to 0$ as $\lambda\to 0$ and such 
that the following holds. If $M'\subset \bC^{N'}$ is a
real-analytic generic submanifold of codimension $d'$ through $0$ and $\2z$ a
germ at $0$ of a holomorphic submersion $\2z\colon (\bC^{N'},0)\to
(\Sigma'_0,0)$, where
$\Sigma'_0$ is the Segre variety of $M'$ at $0$, such that $\2z^{-1}(0)$ is transversal
to $\Sigma'_0$, then there exists a germ at $0$ of a holomorphic mapping
$\Phi\colon ((\Sigma'_0)^{2(d+1)},0)\to (\bC^{N'},0)$ satisfying the following.
If $H\colon (\bC^N,0)\to (\bC^{N'},0)$ is a germ at $0$ of a holomorphic
mapping such that
$H(M)\subset M'$, then
\begin{equation}\Label{e:Hdet}
H=\Phi\circ (\overline{\2z\circ H}\circ h^\lambda_1,\2z\circ H\circ
g^\lambda_1,\ldots, \overline{\2z\circ H}\circ h^\lambda_{d}, \2z\circ H\circ
g^\lambda_d,\overline{\2z\circ H}\circ h^\lambda_{d+1},\2z\circ H ),
\end{equation}
for all $\lambda$ sufficiently small.
\end{Thm}

\begin{Rem}\Label{r:form} {\rm It follows from the proof of Theorem \ref{T:Hdet1} that
there exists an integer $l\geq 0$  such that each of the functions
$g^\lambda_1,\ldots,g^\lambda_{d}$, $h^\lambda_1,\ldots, h^\lambda_{d+1}$ is
given by a convergent power series of the form
\begin{equation}
a_0(\lambda)+\sum_{\alpha\in
\bZ^N_+\setminus\{0\}}a_\alpha(\lambda)\frac{Z^\alpha}{\lambda^{l|\alpha|}},
\end{equation}
where the coefficients $a_0(\lambda),a_\alpha(\lambda)$ are holomorphic in the
unit disk
$\mathbb D$, $a_0(0) = 0$. Another way of expressing this is saying that each of the 
functions $g^\lambda_1,\ldots,g^\lambda_{d}$, $h^\lambda_1,\ldots, h^\lambda_{d+1}$
is given by 
\begin{equation}\Label{e:hat}
g_j^\lambda(Z)=\hat g_j\big(\frac{Z}{\lambda^l},\lambda\big),\quad
h_j^\lambda(Z)=\hat h_j\big(\frac{Z}{\lambda^l},\lambda\big)
\end{equation}
where the $\hat g_1,\ldots,\hat g_{d}$, $\hat h_1,\ldots,
\hat h_{d+1}$ are germs at 0 of holomorphic functions $(\bC^N\times
\bC,0)\to (\bC^N,0)$. }
\end{Rem}

\begin{Rem}{\rm Suppose that $P\colon (\bC^{N'},0)\to (\Sigma'_0,0)$ is a germ at
$0$ of a holomorphic
projection, i.e.\ $P|_{\Sigma'_0}$ is the identity on $\Sigma'_0$. The reader can easily
verify that $\2z:=P$ is a holomorphic submersion satisfying the assumptions of 
Theorem \ref{T:Hdet1}. Conversely, any holomorphic submersion $\2z$ as in
Theorem \ref{T:Hdet1} is a projection up to a local biholomorphism of
$\Sigma'_0$ at $0$.  
 }
\end{Rem}

An immediate corollary of Theorem \ref{T:Hdet1} is the following.

\begin{Cor} \Label{C:Hdet1} Let $M$, $M'$, and $\2z$ be as in Theorem
$\ref{T:Hdet1}$. If $H\colon
(\bC^N,0)\to (\bC^{N'},0)$ is a germ at $0$ of a holomorphic mapping such that
$H(M)\subset M'$, then $H$ is uniquely determined by $\2z\circ H$. 
\end{Cor}

An algebraic reformulation of Corollary \ref{C:Hdet1} can be given as follows. For any
complex manifold $X$ and $p\in X$, let $\mathcal O_X(p)$ denote the ring of germs at
$p$ of holomorphic functions on $X$. For $X=\bC^N$, we write $\mathcal O_N(p)$
instead of $\mathcal O_{\bC^N}(p)$. Recall that if $Y\subset  X$ is a complex analytic
subvariety through $p$, then the ring $\mathcal O_Y(p)$ of germs at $p$ of holomorphic
functions on
$Y$ is given by $\mathcal O_{X}(p)/I(Y)$, where $I(Y)$ denotes the ideal of germs
vanishing on $Y$. Let $H$ be a germ at $p$ of a holomorphic mapping $(X,p)\to
(W,q)$, where $X$ and $W$ are complex manifolds. The mapping $H$ induces a ring
homomorphism
$\Phi_H\colon\mathcal O_W(q)\to \mathcal O_X(p)$, given by $\Phi_H(f)=f\circ H$ for
$f\in \mathcal O_W(q)$. The reader can verify that the following result is a
reformulation of Corollary \ref{C:Hdet1}.

\begin{Thm} \Label{T:Hdet2} Let $M$, $M'$, and $\Sigma'_0$ be as in Theorem
$\ref{T:Hdet1}$ and  denote by $\pi$ the canonical homomorphism $\pi\colon \mathcal
O_{N'}(0)\to \mathcal O_{\Sigma'_0}(0)$. Let  $\phi\colon \mathcal
O_{\Sigma'_0}(0)\to \mathcal O_{N'}(0)$ be any ring homomorphism such that
$\pi\circ\phi\colon \mathcal O_{\Sigma'_0}(0)\to
\mathcal O_{\Sigma'_0}(0)$ is an isomorphism. Then, for any 
$H\colon (\bC^N,0)\to (\bC^{N'},0)$ a germ at $0$ of a holomorphic mapping such that
$H(M)\subset M'$,  the induced homomorphism $\Phi_H$ is uniquely determined by
$\Phi_H\circ \phi$.
\end{Thm}

We give now a necessary and sufficient
condition for a germ at
$0$ of a holomorphic mapping
$F\colon (\bC^N,0)\to (\Sigma'_0,0)$  to be of the form $\2z\circ H$, for some
$\2z$ as in Theorem \ref{T:Hdet1} and a holomorphic mapping $H$ sending 
$M$ into $M'$.

\begin{Thm}\Label{T:F} Let $M$ and $d$ be as in Theorem
$\ref{T:Hdet1}$. Then there exist an open, connected subset $\Omega\subset
\bC^N\times\bC^{2(d+1)(N-d)-N}$ such that the set $$\{(Z,\xi)\in
\bC^N\times\bC^{2(d+1)(N-d)-N}\colon (Z,\xi)\in \Omega,\ Z=0\},$$ is open in
$\{0\}\times\bC^{2(d+1)(N-d)-N}$ with $0$ in its closure, and $2d+1$ germs at
$0$ of holomorphic mappings
$$A_1,\ldots,A_{d},B_1,\ldots, B_{d+1}\colon
(\Omega,0)\to \bC^N$$ such that
$A_j(0,\xi)\to 0$, $B_j(0,\xi)\to 0$ as $\xi\to 0$ (for $(0,\xi)\in \Omega$) and such 
that the following holds. If $M'$, $d'$, $\Sigma_0'$ are as in Theorem
$\ref{T:Hdet1}$, then there exists a germ at $0$ of a holomorphic mapping $\Psi\colon
((\Sigma'_0)^{2(d+1)},0)\to (\bC^{d'},0)$ satisfying the following. Let $F\colon
(\bC^N,0)\to (\Sigma_0',0)$ be a germ at $0$ of a holomorphic mapping. If 
there exists
a germ at $0$ of a holomorphic mapping $H\colon
(\bC^N,0)\to (\bC^{N'},0)$ 
and a germ
at $0$ of a holomorphic submersion $\2z\colon (\bC^{N'},0)\to (\Sigma'_0,0)$ with
$\2z^{-1}(0)$ transversal to $\Sigma'_0$ such that
\begin{equation}\Label{e:HF}
H(M)\subset M',\quad F=\2z\circ H
\end{equation}
then
\begin{multline}\Label{e:Hdet2}
\frac{\partial}{\partial\xi}\Psi (\overline{F}\circ B_1(Z,\xi),F\circ A_1(Z,\xi),\ldots,
 \overline{F}\circ B_{d}(Z,\xi),F\circ A_{d}(Z,\xi),\\ \overline{F}\circ
B_{d+1}(Z,\xi),F(Z) )\equiv 0.
\end{multline}
Conversely, suppose \eqref{e:Hdet2} holds. Then for every germ
at $0$ of a holomorphic submersion $\2z\colon (\bC^{N'},0)\to
(\Sigma'_0,0)$ with
$\2z^{-1}(0)$ transversal to $\Sigma'_0$, there is a unique germ at $0$ of a 
holomorphic mapping $H\colon
(\bC^N,0)\to (\bC^{N'},0)$ satisfying
\eqref{e:HF}.
\end{Thm}

By combining Theorem \ref{T:F} with Theorem 3.1 in \cite{er1}, we obtain the
following result concerning the biholomorphic equivalence problem.

\begin{Cor}\Label{C:F} Let $M$ and $d$ be as in Theorem
$\ref{T:Hdet1}$. Then there exist an open set $\Omega$ and
$A_1,\ldots,A_{d},B_1,\ldots, B_{d+1}$  as in Theorem $\ref{T:F}$ such 
that the following holds. If $M'$ is a real-analytic generic submanifold of codimension
$d$ through $0$ in $\bC^N$, then there exists a germ at
$0$ of a holomorphic mapping $\Psi\colon ((\Sigma'_0)^{2(d+1)},0)\to (\bC^{d},0)$
such that $M$ and $M'$ are biholomorphically equivalent at $0$ if and only if there
exists a germ at $0$ of a holomorphic
mapping $F\colon (\bC^N,0)\to (\Sigma_0',0)$  such that $F|_{\Sigma_0}$ is a local
biholomorphism at $0$ and
\eqref{e:Hdet2} holds.  Here $\Sigma_0$ denotes the Segre variety of $M$ at $0$.
\end{Cor}

We would like to point our that the hypothesis of finite type in Theorems \ref{T:Hdet1}
and \ref{T:F} is crucial as is illustrated by the following simple example.

\begin{Ex} {\rm Let $M\subset \bC^2$ be the real-analytic hypersurface given by 
$$
\im w=(\re w)|z|^2,
$$
which is of finite type at all points except along $\{w=0\}$. Note that the family of
holomorphic mappings
$$
H_t(z,w)=(z,tw),
$$
for all $t\in \bR$, sends $M$ into itself. Thus, the conclusion of Theorem \ref{T:Hdet1}
(with $\2z(z,w)=(z,0)\in \Sigma_0$) does not hold. Also, for any holomorphic
function
$F(z,w)$, with
$F(0)=0$, the mapping
$$
H(z,w)=(F(z,w),0)
$$
sends $M$ into itself and, hence, in contrast with the conclusion of Theorem \ref{T:F},
there is no (non-trivial) condition on a mapping $F$ to be a component of a holomorphic
mapping
$M$ into itself.  }
\end{Ex}

As an application of Theorem \ref{T:Hdet1}, we give a refinement of some results 
concerning finite jet determination of holomorphic mappings between generic
submanifolds (see Section \ref{s:jet} for details). This is a problem that has 
received much attention recently. We mention here the papers
\cite{zautom}, \cite{berconv},
\cite{bmr}, \cite{elz}, \cite{kz1}, \cite{lm2}, where results on finite jet
determination of mappings between generic submanifolds are obtained. The
reader is also referred to the survey papers \cite{z02} and \cite{b} for
further references and results.

\section{An invariant description of normal coordinates}

 For the proof of Theorem \ref{T:Hdet1}, we shall need
the following description of all normal coordinates for a real-analytic generic submanifold.
Let
$M\subset \bC^N$ be a real-analytic generic submanifold through $0$ of codimension
$d$. Recall that local holomorphic coordinates $(z,w)\in \bC^n\times\bC^d$, with
$n=N-d$, are called {\it normal} if $M$ can be expressed near $0$ as a graph of the form
\begin{equation}\Label{e:phi1}
\im w=\phi(z,\bar z,\re w),
\end{equation}
where $\phi(z,\bar z,s)$ is an $\bR^d$-valued real-analytic function in a
neighborhood of $0$ in $\bC^n\times\bR^d$ with
\begin{equation}\Label{e:phi2}
\phi(z,0,s)\equiv\phi(0,\chi,s)\equiv0.
\end{equation}
Equivalently, $M$ can be defined by a complex equation of the form 
\begin{equation}\Label{e:Q}
w=Q(z,\bar z,\bar w),
\end{equation}
where $Q(z,\chi,\tau)$ is a $\bC^d$-valued holomorphic function, defined near $0$ in
$\bC^n\times\bC^n\times \bC^d$, satisfying
\begin{equation}\Label{e:Qnormal}
Q(z,0,\tau)\equiv Q(0,\chi,\tau)\equiv \tau.
\end{equation}
Normal coordinates were first introduced by Chern--Moser in \cite{cm} (see also
\cite{book}).

As in the beginning of Section
\ref{s:intro}, let
$U$ be a small open neighborhood of $0$ in $\bC^N$, and  $\Sigma_0$ the Segre variety
of $M$ at
$0$. Let
$\2z\colon U\to
\Sigma_0$ be a holomorphic submersion such that
$\2z(0)=0$ and the $d$-dimensional complex submanifold $W:=\2z^{-1}(0)$ is
transversal to
$\Sigma_0$ at $0$.  Observe that, for $p\in \bC^N$ sufficiently close to $0$, the
submanifolds $\Sigma_p$ and $W$ also intersect transversally near $0$.
Hence, after shrinking
$U$ if necessary, we
may define a mapping
$\sigma\colon U\to W$ by letting $\sigma(p)$ be the single point of intersection between
$\Sigma_p$ and $W$.  We denote by $\iota\colon W\to W$ the restriction of $\sigma$ to
$W$, i.e.
\begin{equation}\Label{e:iota}
\Sigma_p\cap W=\{\iota(p)\},\quad p\in W.
\end{equation} 
It follows that, for $p\in W$ sufficiently close to $0$, $\iota^2(p)=p$ since $\iota(p)\in
\Sigma_p\cap W$ and, hence, $\iota $ is a local involution on $W$.

\begin{Thm}\Label{t:normal1}  Let $M$ be a real-analytic generic submanifold of
codimension $d$ in $\bC^N$ with $0\in M$, and $\Sigma_p$ its Segre variety at $p$,
for $p$ close to $0$.  Let $U\subset \bC^N$ be a sufficiently small neighborhood of $0$ 
and $\2z\colon (U,0)
\to (\Sigma_0,0)$ be a holomorphic submersion such that the $d$-dimensional complex
submanifold
$W:=\2z^{-1}(0)$ is transversal to $\Sigma_0$ at $0$, and $\iota\colon W\to W$ the
corresponding mapping as defined in \eqref{e:iota}.  Then there are open neighborhoods
of the origin $V\subset \bC^N$, $X\subset \Sigma_0$, $Y\subset W$ such that the
following hold:
\begin{enumerate}
\item [(i)] $\iota\colon Y\to Y$ is an anti-holomorphic involution fixing $M\cap Y$.

\item [(ii)] There is a unique holomorphic submersion $\2w\colon (V,0)\to (Y,0)$ such
that the mapping $H\colon (V,0)\to (X\times Y,0)$, where $H(Z)\:=(\2z(Z),\2w(Z))$, is a biholomorphism satisfying
\begin{equation}\Label{e:Hsigma}
H(\Sigma_p\cap V)=X\times\{\iota(p)\},\quad \forall p\in Y.
\end{equation}

\item[(iii)] If $\alpha\colon (X,0)\to (\bC^n,0)$, and $\beta\colon (Y,0)\to (\bC^d,0)$ are
biholomorphisms, then in the coordinates $(z,w):=(\alpha\circ \2z,\beta\circ\2w)$ in $\bC^N=\bC^n\times\bC^d$ the submanifold $M$ is given near $0$ by 
$
w=Q(z,\bar z,\bar w),
$
where $Q(z,\chi,\tau)$ is a $\bC^d$-valued holomorphic function satisfying
\begin{equation}\Label{e:Qsubnormal}
Q(z,0,\tau)\equiv Q(0,\chi,\tau)\equiv \hat\iota(\bar \tau)
\end{equation}
with $\hat\iota$  the involution given by $\beta\circ\iota\circ\beta^{-1}$. Moreover, if
$\beta(Y\cap M)\subset \bR^d$, then $\hat \iota(w)=\bar w$ and, hence, $(z,w)$ are
normal coordinates, i.e.\ the identity
\eqref{e:Qnormal} holds.
\end{enumerate}
\end{Thm}

\begin{proof}[Proof of Theorem $\ref{t:normal1}$]
We let $\rho(Z,\bar Z)$ be a defining
function for $M$ as in the introduction. Consider the germ at $0$ of a holomorphic
mapping $f\colon (\bC^N\times
\bC^N,0)\to (\bC^d\times\Sigma_0,0)$ defined by 
\begin{equation}
f(Z,\zeta):=(\rho(Z,\zeta),\2z(Z)),
\end{equation}
and the equation
\begin{equation}\Label{e:Feq}
f(Z,\zeta)=(0,t).
\end{equation}
We claim that $Z\mapsto f(Z,0)$ is a local biholomorphism at $0$. Indeed, 
\begin{equation}
\frac{\partial f}{\partial Z}(0,0)=\bigg( \frac{\partial \rho}{\partial Z}(0,0), \frac{\partial
\2z}{\partial Z}(0)\bigg)
\end{equation}
and, hence, the claim follows from the transversality of the intersection between
$\Sigma_0$ and $W$ at $0$. By the implicit function theorem, there exists a unique germ
at $0$ of a holomorphic mapping $Z=\gamma(t,\zeta)$ from
$(\Sigma_0\times\bC^N,0)\to (\bC^N,0)$  that solves equation \eqref{e:Feq}.  It follows
from \eqref{e:Feq} that $t\mapsto \gamma(t,0)$ has rank $n:=N-d$ at $0$ and that
$t\mapsto
\gamma(t,\bar p)$, for $p\in\bC^N$ sufficiently close to $0$, parametrizes an open
piece of the Segre variety $\Sigma_p$. We observe, from the definition of $\gamma$, 
that $\sigma(p)=\gamma(0,\bar p)$, 
for $p\in \bC^N$ close to $0$; recall that $\sigma(p)$ denotes the single point of
intersection between $\Sigma_p$ and $W$, as defined above. In particular, 
\begin{equation}
\iota(p)=\gamma(0,\bar p),
\end{equation}
where $\iota$ is the involution of $W$ defined by \eqref{e:iota}, and hence 
the anti-holomorphic mapping $p\mapsto \gamma(0,\bar p)$ is a local diffeomorphism at
$0$ of $W$. It follows that the mapping $(t,p)\mapsto \gamma(t,\bar p)$ from
$\Sigma_0\times W\to \bC^N$ is holomorphic in $t$, anti-holomorphic in $p$, and is a
local diffeomorphism at $0$. Hence, if we denote by $W^*$ the submanifold $\{Z\colon
\bar Z\in W\}$, then the mapping  $\Gamma(t,p):=\gamma(t, p)$ from $\Sigma_0\times
W^*\to \bC^N$ is a local biholomorphism at $0$. As a consequence, we may define the
germ at
$0$ of a holomorphic mapping 
$\2w\colon (\bC^N,0) \to (W,0)$ by 
\begin{equation}\Label{e:wtilde}
\2w(\Gamma(t,p))=\gamma(0,p)\ (=\iota(\bar p)).
\end{equation}
Since $p\mapsto \iota(p)$ is a local diffeomorphism of $W$ at $0$, it follows
that $\2w$ is a submersion at $0$. Since $\iota$ is a local involution on $W$, we can
find a sufficiently small open neighborhood $Y$ of $0$ in $W$ such that $\iota$ is an
involution on $Y$ (i.e.\ $\iota$ maps $Y$ onto itself and $\iota^2$ is the identity).
Indeed, if $Y_0$ is any sufficiently small open neighborhood of $0$ in $W$, then
$Y:=Y_0\cap \iota(Y_0)$ is such a neighborhood. This proves (i).
To prove (ii), let $X$ be any sufficiently small open
neighborhood of $0$ in $\Sigma_0$ and define
$V:=\Gamma(X\times Y^*)$. Then $\2w$, defined by \eqref{e:wtilde} is a
holomorphic submersion
$(V,0)\to (Y,0)$. Observe that, for every $p\in Y$,  $X\ni t\mapsto \Gamma(t,\bar p)\in V$
parametrizes $\Sigma_p\cap V$. Hence, equation \eqref{e:Hsigma} in Theorem
\ref{t:normal1} is equivalent to
\eqref{e:wtilde}.  This proves (ii).

To prove (iii), we assume that $\alpha$ and $\beta$ are as in (iii), and let $(z,w)$ be the
coordinates $(z,w):=(\alpha\circ\2z,\beta\circ\2w)$. In these coordinates, it
follows from (ii) that $\Sigma_0=\{(z,w)\colon w=0\}$. Consequently, if
$\hat\rho(z,w,\bar z,\bar w)=0$ is a defining equation for $M$ in the coordinates $(z,w)$,
then $\det\, \partial\hat \rho/\partial w(0)\neq 0$. Hence, by the implicit function theorem,
we may solve for $w$ in the equation $\hat\rho(z,w,\bar z,\bar w)=0$ and obtain a
defining equation for $M$ of the form \eqref{e:Q}. The fact that 
\begin{equation}
\Label{e:halfQ}Q(z,0,\tau)=\hat \iota
(\bar\tau)
\end{equation}
 is a direct consequence of \eqref{e:Hsigma}. To prove the remaining part
of
\eqref{e:Qsubnormal}, we note, by the fact that \eqref{e:Q} defines a real submanifold,
that we have
\begin{equation}\Label{e:reality2}
Q(z,\chi,\bar Q(\chi,z,w))\equiv w.
\end{equation}
By substituting $z=0$ in \eqref{e:reality2} and using \eqref{e:halfQ}, we obtain
$Q(0,\chi,\overline{\hat \iota(w)})=w$. The desired identity
$Q(0,\chi,\tau)=\hat\iota(\bar\tau)$ follows by taking $w=\hat\iota(\bar\tau)$. 

If $\beta(Y\cap M)\subset \bR^d$, then, since $Y\cap M$ is the fixed point set of the
involution $\iota$, it follows that $\bR^d\cap \beta(Y)$ is the fixed point set of the
anti-holomorphic involution $\hat\iota$ on $\beta(Y)\subset \bC^d$, i.e.\ $\hat
\iota(w)=\bar w$. The identity \eqref{e:Qnormal} follows immediately. The fact that 
$M$ can be graphed as in \eqref{e:phi1} with $\phi$ satisfying \eqref{e:phi2} is a direct
consequence of the implicit function theorem and \eqref{e:Qnormal}. This completes the
proof of Theorem \ref{t:normal1}.
\end{proof}

\begin{Rem}{\rm It is not difficult to see that all normal coordinates $(z,w)$ are obtained
in the way described by Theorem \ref{t:normal1} for some choice of submersion $\2z$. The details of this are left to the reader. }
\end{Rem}

\begin{Rem}{\rm Let $Z=(\hat z,\hat w)$ be given coordinates in $\C^N =
\C^n\times\C^d$ in which the Segre variety $\Sigma_0$ of $M$ at $0$ is tangent to
$\{(\hat z,\hat w)\colon
\hat w=0\}$ at $0$. As a consequence of Theorem \ref{t:normal1} (and its proof), we
obtain the following description of all possible holomorphic transformations
$(z,w)=(F(\hat z,\hat w),G(\hat z,\hat w))$ yielding normal coordinates for $M$. Let
$F\colon (\bC^N,0)\to (\bC^n,0)$ be an arbitrary local holomorphic mapping with
$\det\, (\partial F/\partial\hat z)(0)\neq0$. In the setting of the theorem, this corresponds
to a choice of a holomorphic submersion $\2z$ and a local chart $\alpha$ of
$\Sigma_0$, with $F=\alpha\circ\2z$. We note that $\2z$ and $\alpha$ are not
uniquely determined by $F$. However, any two different choices of $\2z$ differ only
by a reparametrization of $\Sigma_0$. An inspection of the proof of Theorem
\ref{t:normal1} shows that the $\2w$ given by (ii) in the theorem is uniquely
determined by
$F$.  Moreover, for any local chart
$\beta$ on $W$ such that $\beta(M\cap W)\subset \bR^d$, the mapping $G=\beta\circ
\2w$ produces normal coordinates by (iii) of Theorem \ref{t:normal1}.  It is easily
seen that if we write
$G(0,\hat w)=g_1(\hat w)+ig_2(\hat w)$, where $g_1$ and $g_2$ are real-valued on
$\bR^d$, then there is a one-to-one correspondence between choices of such
parametrizations $\beta$ and choices of $g_1(\hat
w)$ with $\det\,(\partial g_1/\partial \hat w)(0)\neq 0$. We conclude that 
$G(\hat z,\hat w)$ is uniquely determined by $F(\hat z,\hat w)$ and an  arbitrary choice of
$g_1(\hat w)$ with $\det\,(\partial g_1/\partial \hat w)(0)\neq 0$. 
}
\end{Rem}

\section{Proof of Theorem $\ref{T:F}$ in the case of hypersurfaces
in $\bC^2$}
\Label{s:hyper}

In order to illustrate the idea of the proof of Theorem \ref{T:F}, we first give a proof for
the case of hypersurfaces $M,  M'$ in $\bC^2$. Let 
$M\subset
\C^2$ be a real-analytic hypersurface of finite type at
$0\in M$.  Assume that
$(z, w)\in \C\times\C$ are normal coordinates at 0. Thus $M$ is given locally near $0$ 
by \eqref{e:Q}, where the
scalar-valued holomorphic function
$Q(z,\chi,\tau)$ satisfies \eqref{e:Qnormal}. The finite type condition on the
hypersurface $M$ is equivalent to $Q(z,\chi,0)\not\equiv 0$, which implies, by the
normality of $(z,w)$ that
\begin{equation}\Label{e:ft}
Q_\chi(z,\chi,0)\not\equiv 0,
\end{equation}
where we use the notation $Q_\chi=\partial Q/\partial \chi$. The first four iterated  Segre
mappings  (as defined in \cite{BERrat} and \cite{BERJAG}) are given by:
\begin{equation}
\begin{aligned}
v^1(t^1) &:=(t^1,0)\\
v^2(t^1,t^2)&:=(t^2,Q(t^2,t^1,0))\\
v^3(t^1,t^2,t^3)&:=(t^3,Q(t^3,t^2,\bar Q(t^2,t^1,0)))\\
v^4(t^1,t^2,t^3,t^4)&:=(t^4,Q(t^4,t^3,\bar Q(t^3,t^2,Q(t^2,t^1,0)))).
\end{aligned} 
\end{equation}
For convenience, we shall also write $v^k(t^1,\ldots, t^k)=(t^k,u^k(t^1,\ldots, t^k))$. 

We let $M'\subset \bC^2$ be another real-analytic hypersurface through $0$
and $(\2z,\2w)$ normal coordinates for $M'$. We refer to the
corresponding objects for $M'$ by the addition of $\ \widetilde{}$ . Let
$H=(F,G)$ be a germ at $0$ of a holomorphic mapping with $H(0)=0$. If $H$ sends
$M$ into $M'$, then (see \cite{er1}, Section 2)
\begin{equation}\Label{e:Gv4}
G\circ v^4=\2u^4(\overline{F\circ v^1}, F\circ v^2, \overline{F\circ v^3}, F\circ v^4).
\end{equation}
Conversely, if $H$ satisfies \eqref{e:Gv4}, then we claim that $H$ sends $M$ into
$M'$. Indeed, if we take $t^1=0$ in \eqref{e:Gv4} and complex conjugate, 
then we obtain
\begin{equation}\Label{e:Gv3}
\overline{G\circ v^3}=\overline{\2u^3}(F\circ v^1, \overline{F\circ v^2}, F\circ v^3)
\end{equation}
by using standard properties of the iterated Segre mappings (see \cite{BERrat} or
\cite{BERJAG}). We now observe that
\begin{equation}\Label{e:newid}
\2u^4(\tilde t^1,\tilde t^2,\tilde t^3, \tilde t^4)=\widetilde Q(\tilde t^4,\tilde
t^3,\overline{\2u^3}(\tilde t^1,\tilde t^2,\tilde t^3)).
\end{equation}
By  using \eqref{e:newid} and \eqref{e:Gv3} in \eqref{e:Gv4}, we conclude that 
\begin{equation}\Label{e:Gv4Q}
G\circ v^4=\widetilde Q(F\circ v^4, \overline{F\circ v^3}, \overline{G\circ
v^3}).
\end{equation}
Let $\mathcal M\subset \bC^2\times\bC^2$ be the complexification of $M$, i.e.\ the
complex submanifold through $0$ in $\bC^2\times\bC^2$ defined by
$w=Q(z,\chi,\tau)$. Since $(t^1,t^2,t^3,t^4)\mapsto
(v^4(t^1,t^2,t^3,t^4),\overline{v^3}(t^1,t^2,t^3))$ is a holomorphic mapping of generic
full rank  into $\mathcal
M$ (see \cite{BERrat} or \cite{BERJAG}), we
conclude that 
$G(z,w)=\widetilde Q(F(z,w),\bar F(\chi,\tau),\bar G(\chi,\tau))$ for all $(z,w,\chi,\tau)$
on
$\mathcal M$ and, hence,
$H$ sends
$M$ into $M'$. This proves the claim.  

Consider the equation $(z,w)=v^4(t^1,t^2,t^3,t^4)$, which can also be written in the
form $z=t^4$ and 
\begin{equation}
w=Q(z,t^3,\bar Q(t^3,t^2,Q(t^2,t^1,0))).
\end{equation}
We make the linear change of variables
\begin{equation}\Label{e:linchange}
\eta^1=\frac{t^1+t^3}{2},\quad \eta^2=t^2,\quad \sigma=\frac{t^1-t^3}{2}
\end{equation}
and obtain 
\begin{equation}\Label{e:wisQ}
w=Q(z,\eta^1-\sigma,\bar Q(\eta^1-\sigma,\eta^2,Q(\eta^2,\eta^1+\sigma,0))).
\end{equation}
Let us use the notation $\eta=(\eta^1,\eta^2)$ and write
$$U(\eta,z,\sigma):=Q(z,\eta^1-\sigma,\bar
Q(\eta^1-\sigma,\eta^2,Q(\eta^2,\eta^1+\sigma,0))).
$$ We have
$U(\eta,0,0)\equiv 0$ and 
\begin{multline}\Label{e:210}
\Delta(\eta):=\frac{\partial}{\partial \sigma}U(\eta,z,\sigma)\big|_{z=\sigma=0}=\\
-\bar
Q_\chi(\eta^1,\eta^2,Q(\eta^2,\eta^1,0))+\bar
Q_w(\eta^1,\eta^2,Q(\eta^2,\eta^1,0))Q_\chi(\eta^2,\eta^1,0).
\end{multline}
Here we have used the notation $\bar Q(\chi,z,w)$ and hence the corresponding
derivatives $\bar Q_\chi$  and $\bar Q_w$ refer to partial derivatives with respect to the
first and last variable, respectively.  By
differentiating the identity
$$
\bar Q(\eta^1,\eta^2,Q(\eta^2,\eta^1,0))\equiv 0
$$
with respect to $\eta^1$, we obtain from \eqref{e:210}
$$
\Delta(\eta)=2\bar Q_w(\eta^1,\eta^2,Q(\eta^2,\eta^1,0))Q_\chi(\eta^2,\eta^1,0)\not\equiv
0,
$$
in view of  \eqref{e:ft} and the fact that $\bar Q_w(0,0,0)=1$. We may now apply
the singular implicit function theorem given in Proposition 4.1.18 in \cite{BERrat} and
conclude that the equation \eqref{e:wisQ} has a unique solution of the form
\begin{equation}\Label{e:sigma}
\sigma=\Theta\big(\eta,\frac{z}{\Delta(\eta)^2}, \frac{w}{\Delta(\eta)^2}\big ),
\end{equation}
where $\Theta(\eta,z',w')$ is holomorphic near $0\in \bC^4$ and
$\Theta(\eta,0,0)\equiv0$. If we now substitute for $(t^1,t^2,t^3,t^4)$ in \eqref{e:Gv4}
using the linear change of variables \eqref{e:linchange}, $t^4=z$, and then substitute for
$\sigma$ using
\eqref{e:sigma}, then we obtain
\begin{multline}\Label{e:GuF}
G(z,w)=\\\2u^4(\overline{F\circ v^1}(\eta^1+\Theta), F\circ
v^2(\eta^1+\Theta,\eta^2),
\overline{F\circ v^3}(\eta^1+\Theta,\eta^2,\eta^1-\Theta), F(z,w)),
\end{multline}
where $\Theta$ is given by the right hand side of \eqref{e:sigma}. In particular, if $H$
sends
$M$ into
$M'$, then the right hand side of  \eqref{e:GuF} is independent of the
variable $\eta$. Conversely, if $F(z,w)$ is such that the right hand side of  \eqref{e:GuF}
is independent of the variable $\eta$, then we can define $G(z,w)$ by \eqref{e:GuF}. We
claim that $H=(F,G)$ sends $M$ into $M'$. Indeed, for any $\eta\in \bC^2$
sufficiently close to $0$ with $\Delta(\eta)\neq 0$, we have 
\begin{equation}\Label{e:identity}
\Theta\big (\eta,\frac{z}{\Delta(\eta)^2},\frac{U(\eta,z,\sigma)}{\Delta(\eta)^2}\big
)=\sigma,
\end{equation}
for all sufficiently small $z$ and $\sigma$, 
by the uniqueness of the solution \eqref{e:sigma} to the equation \eqref{e:wisQ}. 
We now make the subsitution 
$(z,w)=v^4(t^1,t^2,t^3,t^4)$, $\eta^1=(t^1+t^3)/2$, $\eta^2=t^2$ in
\eqref{e:GuF}. Using again the linear change of variables \eqref{e:linchange}, and 
$t^4=z$
in the identity \eqref{e:identity}, we conclude, since in these variables we have
$u^4(t^1,t^2,t^3,t^4)=U(\eta,z,\sigma)$, that
\eqref{e:Gv4} holds. This proves the claim, in view of the remarks above. 

This proves Theorem \ref{T:F} for hypersurfaces in $\bC^2$ with 
\begin{multline}
\Psi (\overline{F}\circ B_1(Z,\eta), F\circ A_1(Z,\eta), 
  \overline{F}\circ B_{2}(Z,\eta), F(Z))=\\
\2u^4(\overline{F\circ v^1}(\eta^1+\Theta), F\circ
v^2(\eta^1+\Theta,\eta^2),
\overline{F\circ v^3}(\eta^1+\Theta,\eta^2,\eta^1-\Theta), F(z,w)),
\end{multline}
and hence
$$
A_1(z,w,\eta)=v^2(\eta^1+\Theta,\eta^2),\ B_1(z,w,\eta)=\overline{v^1}(\eta^1+\Theta),\
B_2(z,w,\eta)=\overline{v^3}(\eta^1+\Theta,\eta^2,\eta^1-\Theta),
$$
where $\Theta$ is given by the right hand side of \eqref{e:sigma}. The open set
$\Omega\subset \bC^2\times\bC^2$ can be taken to be of the form
$$
 |\eta|<\epsilon,\quad \Delta(\eta)\neq 0,\quad |z|<\epsilon \Delta(\eta),\quad |w|<\epsilon
\Delta(\eta),
$$
for $\epsilon>0$ sufficiently small.

\section{Proof of Theorem $\ref{T:F}$ in the general case and the proof of Theorem~
$\ref{T:Hdet1}$}

\begin{proof}[Proof of Theorem $\ref{T:F}$] We point out that, in view of Theorem
\ref{t:normal1}, it suffices to prove Theorem \ref{T:F} in some sets of  normal
coordinates for
$M$ and
$M'$.  The proof in the general case parallels that
for hypersurfaces in $\bC^2$ given in Section \ref{s:hyper}. Let $M\subset
\bC^N$  be given in normal coordinates
$(z,w)\in \bC^n\times\bC^d$  near the
origin in $\bC^N$, i.e.\ by \eqref{e:Q} where $Q(z,\chi,\tau)$ is a $\bC^d$-valued
holomorphic function satisfying \eqref{e:Qnormal}. The iterated Segre mappings (see
\cite{BERrat} and \cite{BERJAG}) $v^j\colon (\bC^{jn},0)\to (\bC^N,0)$ are given in
these coordinates by $v^1(t^1)=(t^1,0)$ and, recursively, for $j\geq 1$ by 
\begin{equation}
v^{j+1}(t^1,\ldots, t^{j+1})=(t^{j+1}, u^{j+1}(t^1,\ldots,
t^{j+1}))=(t^{j+1},Q(t^{j+1},\overline{v^j}(t^1,\ldots, t^j)).
\end{equation}

We let $M'\subset \bC^{N'}$ be another real-analytic generic submanifold
through
$0$ and $(\2z,\2w)$ normal coordinates for $M'$. We refer to the
corresponding objects for $M'$ by the addition of $\ \widetilde{}$ . Let
$H=(F,G)$ be a germ at $0$ of a holomorphic mapping $(\bC^N,0)\to (\bC^{N'},0)$ 
with $H(0)=0$. Let $m:=d+1$. As in Section \ref{s:hyper}, if $H$ sends
$M$ into $M'$ (cf.\ \cite{er1}, Section 2), then
\begin{equation}\Label{e:Gv2m}
G\circ v^{2m}=\2u^{2m}(\overline{F\circ v^1}, F\circ v^2, \ldots,\overline{F\circ
v^{2m-1}}, F\circ v^{2m}).
\end{equation}
Conversely, if \eqref{e:Gv2m} holds, then, since the mapping 
$$
(t^1,\ldots, t^{2m})\mapsto (v^{2m}(t^1,\ldots, t^{2m}), \overline{v^{2m-1}}(t^1,\ldots,
t^{2m-1}))
$$
has generic full rank as a holomorphic mapping into the complexification $\mathcal
M\subset
\bC^N\times
\bC^N$ (see \cite{BERrat} or \cite{BERJAG}), a similar argument to that in Section
\ref{s:hyper} shows that $H$ sends
$M$ into $M'$. 

As in Section \ref{s:hyper}, we consider the equation $(z,w)=v^{2m}(t^1,\ldots,
t^{2m})$ or, equivalently, $z=t^{2m}$ and 
\begin{equation}\Label{e:wisu}
w=u^{2m}(t^1,\ldots,t^{2m-1},z).
\end{equation}
We make the linear change of variables 
\begin{equation}\Label{e:tlinear}
\eta^j=\frac{t^j+t^{2m-j}}{2},\ \sigma^j=\frac{t^j-t^{2m-j}}{2},
\quad j=1,\ldots, m-1,
\end{equation}
and
\begin{equation}
\eta^m=t^m.
\end{equation}
Thus, if we write $\eta=(\eta^1,\ldots,\eta^m)$, $\sigma=(\sigma^1,\ldots,
\sigma^m)$, and 
\begin{equation}
U(\eta,z,\sigma):=u^{2m}(\eta^1+\sigma_1,\ldots, \eta^{m-1}+\sigma_{m-1},\eta^m,
\eta^1-\sigma_1,\ldots, \eta^{m-1}-\sigma_{m-1},z)
\end{equation}
then equation \eqref{e:wisu} becomes
\begin{equation}
w=U(\eta,z,\sigma).
\end{equation}
Since $M$ is of finite type at $0$, it follows from Lemma 4.1.3 in \cite{BERrat} that
$U(\eta,0,0)\equiv 0$ and that we may decompose $\sigma=(\sigma',\sigma'')\in
\bC^d\times\bC^{(m-1)n-d}$, after reordering the variables if necessary, such that
\begin{equation}\Label{e:delta}
\Delta(\eta):=\det
\big(\frac{\partial}{\partial\sigma'}U(\eta,z,\sigma',\sigma'')\big|_{z=0,\sigma=0}\big)
\not\equiv 0.
\end{equation}
(The reader
should be warned that the notation in \cite{BERrat} is slightly different from that of the
present paper.)  We should point out that in the hypersurface case in Section \ref{s:hyper}
these facts were easily verified directly by using the condition of finite type. In the case 
of higher codimension, the proof of these facts is more involved. We may now apply the
singular implicit function theorem given in Proposition 4.1.18 of \cite{BERrat} and solve
for $\sigma'$ in equation \eqref{e:wisu} and obtain a unique solution of the form
\begin{equation}\Label{e:sigma'}
\sigma'=\Theta(\eta, \frac{\sigma''}{\Delta(\eta)^2},\frac{z}{\Delta(\eta)^2},
\frac{w}{\Delta(\eta)^2}),
\end{equation}
where $\Theta(\eta,\sigma'',z,w)$ is a $\bC^d$-valued  holomorphic function near $0$ in
$\bC^{mn}\times\bC^{(m-1)n-d}\times\bC^{n}\times\bC^d$ with
$\Theta(\eta,0,0,0)=0$. Substituting for $(t^1,\ldots, t^{2m})$ in terms of $\eta,$
$\sigma$, and $z$ in
\eqref{e:Gv2m} using
\eqref{e:tlinear}, $t^m=\eta^m$, $t^{2m}=z$, and then substituting for $\sigma'$ using
\eqref{e:sigma'}, we obtain a relation of the form (cf.\  \eqref{e:GuF})
\begin{multline}\Label{e:GuF2}
G(z,w)=\Psi (\overline{F}\circ B_1(z,w,\xi),F\circ A_1(z,w,\xi),\ldots,\\
 \overline{F}\circ B_{d}(z,w,\xi),F\circ A_{d}(z,w,\xi), \overline{F}\circ
B_{d+1}(z,w,\xi),F(z,w) )
\end{multline}
where $\xi=(\eta,\sigma'')\in \bC^{nm}\times\bC^{(m-1)n-d}$. The $A_i$ and $B_j$ are
obtained by making the substitutions described above for $t^1,\ldots, t^{2m}$ in the
iterated Segre mappings $v^1,\ldots,v^{2m}$ and their complex conjugates. The reader
can verify from the construction of $\Theta$ that $A_i(0,\xi)$ and $B_j(0,\xi)$ tend to 0
as $\xi\to 0$. If
$H=(F,G)$ sends $M$ into $M'$, then the right hand side of \eqref{e:GuF2} is
independent of $\xi$ or, equivalently, \eqref{e:Hdet2} holds. The converse follows in the
same way as in the hypersurface case in Section \ref{s:hyper}. This completes the proof
of Theorem \ref{T:F}. 
\end{proof}

\begin{proof}[Proof of Theorem $\ref{T:Hdet1}$] We keep the notation of the proof of
Theorem \ref{T:F}. First, by using Theorem \ref{t:normal1}, we can find normal
coordinates $(\2z,\2w)$ for $M'$ such that $\2z=\2z(\2z)$ is the
submersion given in Theorem \ref{T:Hdet1}. If $H$ is a germ at
$0$ of a holomorphic mapping $(\bC^N,0)\to (\bC^{N'},0)$, then, in  these coordinates,
$H=(F,G)$ with $F=\2z\circ H$.   Hence, to prove Theorem~\ref{T:Hdet1}, it
suffices to show that we have an identity of the form
\begin{equation}\Label{e:Gdet}
G=\Psi\circ (\overline{F}\circ h^\lambda_1,F\circ
g^\lambda_1,\ldots, \overline{F}\circ h^\lambda_{d}, F\circ
g^\lambda_d,\overline{F}\circ h^\lambda_{d+1},F ),
\end{equation}
for all $\lambda$ sufficiently small.
Let
$D\colon
\mathbb D\to
\bC^{nm}$ be a holomorphic mapping with
$D(0)=0$ such that $\Delta(D(\lambda))\neq 0$ for $\lambda\neq0$, where $\Delta(\eta)$
is the determinant given by \eqref{e:delta}. The conclusion of Theorem \ref{T:Hdet1} now
follows by substituting $\xi=(\eta,\sigma'')=(D(\lambda),0)$ in the identity
\eqref{e:GuF2}. \end{proof}

\begin{Rem}\Label {r:Levi}{\rm   In the case where $M$ is a Levi-nondegenerate
 hypersurface it is possible to prove a version of Theorem \ref{T:F} with fewer
parameters by using the iterated Segre map $v^3$ rather than $v^4$. 
We shall illustrate this in the model case where $M\subset \C^{n+1}$
and $M' \subset \C^{n'+1}$ are respectively given by $$\im w =
<z,\1z>_n =
\sum_{j=1}^n |z_j|^2\quad \text{and} \quad \im \2w =
<\2z,\1{\2z}>_{n'} =
\sum_{j=1}^{n'} |\2z_j|^2.
$$ 
By a calculation similar to that given in Section \ref{s:hyper}, one obtains
that a germ at $0$ of a mapping $H =(F,G): (\C^{n+1},0) \to
(\C^{n'+1},0) $ sends $M$ into
$M'$ if and only if
\begin{multline}\Label{e:ball}
G(z,w) = \\ 2i\big<F(z,w) - F\big(\frac {<z,t^2>_n -
<t^{1}_*,t^{2}_*>_{n-1} + iw/2}{t^2_1},t^{1}_*,0\big),
 \1F(t^2,w-2i<z,t^2>_n)\big >_{n'},
\end{multline}
where we have used the notation $t^j = (t^j_1,t^{j}_*) \in
\C\times\C^{n-1}$, for $j=1,2$.  Here equation \eqref{e:Hdet2} is equivalent to
the condition that the right hand side of \eqref{e:ball} is independent
of the $2n-1$ parameters $t^{1}_*$ and $t^2$. Note that Theorem \ref{T:F}
involves $3n-1$ parameters in the case that $M$ is a hypersurface. 
 }
\end{Rem}

\begin{Rem}{\rm  We point out that the independence of the right hand side
of \eqref{e:Gdet} on the parameter $\lambda$ is not sufficient to guarantee
that $F$ is a component of a mapping sending $M$ into $M'$. This is the
case even in the context of self mappings of the Levi-nondegenerate hypersurface $M$
given in Remark \ref{r:Levi} above
with
$n>1$. }
\end{Rem}

\section{An application to the problem of finite jet determination}\Label{s:jet}

Let $\mathcal F$ be a class of germs of holomorphic mappings $(X,x)\to (Y,y)$, where
$X$ and $Y$ are complex manifolds with $x\in X$ and $y\in Y$, respectively.
We shall say that $\mathcal F$ satisfies the {\it finite jet determination
property} at $x$ if there exists an integer $k\geq 0$ such that for any
pair $H^1,H^2\in\mathcal F$, the condition $j^k_xH=j^k_xH'$ implies
$H\equiv H'$. Here, $j^k_xH$ denotes the $k$-jet at $x$ of $H$. For instance, if $M$
and $M'$ are real-analytic generic submanifolds  of codimension $d$ through $0$ in
$\bC^N$ and $\bC^{N'}$ respectively and $\mathcal F$ is a class of germs at $0$ of
holomorphic mappings
$(\bC^N,0)\to (\bC^{N'},0)$ sending
$M$ into $M'$, then there are a number of sufficient conditions that can be imposed on
$M$ (or
$M'$) to guarantee that $\mathcal F$ satisfies the finite jet determination
property (see e.g. \cite{zautom}, \cite{berconv},
\cite{bmr}, \cite{elz}, \cite{kz1}, \cite{lm2}). As a consequence of Theorem
\ref{T:Hdet1}, we obtain the following
result. 

\begin{Thm}\Label{T:fd} Let $M$ and $M'$ be real-analytic generic submanifolds
through
$0$ in
$\bC^N$ and $\bC^{N'}$, respectively, with $M$ of finite type at $0$. Let $\mathcal F$
be a class of germs at $0$ of holomorphic mappings $(\bC^N,0)\to
(\bC^{N'},0)$ sending
$M$ into $M'$ such that $\mathcal F$ satisfies the finite jet determination
property. Then there exists an integer $k> 0$ with the following property.
Let $\2z$ be a germ
at $0$ of a holomorphic submersion $\2z\colon (\bC^{N'},0)\to (\Sigma'_0,0)$, 
where
$\Sigma'_0$ is the Segre variety of $M'$ at $0$, such that $\2z^{-1}(0)$ is transversal
to $\Sigma'_0$. If $H^1,H^2\in \mathcal F$ and $j^k_0(\2z\circ H^1)=j^k_0(\2z\circ H^2)$, then $H^1\equiv H^2$.
\end{Thm}

By using a result from \cite{bmr}, we immediately obtain the following
corollary. Recall that a generic submanifold is called {\it holomorphically
nondegenerate} at $0$ if there are no germs at $0\in M$ of (non-trivial)
holomorphic vector fields (i.e.\ $(1,0)$ vector fields with holomorphic
coefficients) that are tangent to $M$ in a neighborhood of
$0$. 

\begin{Cor} Let $M$ and $M'$ be real-analytic generic submanifolds of
codimension $d$ through
$0$ in
$\bC^N$ with $M$ of finite type and holomorphically nondegenerate at $0$. 
Let $(\2z,\2w)\in \bC^{N-d}\times\bC^d$ be normal coordinates for $M'$ at $0$.
Then there exists an integer $k\geq 0$ with the following property. Let $H^1, H^2\colon 
(\bC^N,0)\to (\bC^N,0)$ be germs at $0$ of local biholomorphisms sending $M$ into
$M'$ and $H^j=(F^j,G^j)$, $j=1,2$, in the coordinates $(\2z,\2w)$. If $j^k_0
F^1=j^k_0F^2$, then $H^1\equiv H^2$.
\end{Cor}

Recall that a germ at $0$ of a real-analytic hypersurface $M \subset \C^N$ 
is of D'Angelo finite type \cite{d'a} if there is no germ of a nontrivial complex curve
$\mathcal C$ through $0$ contained in $M$.  By using a recent result of Lamel-Mir
\cite{lm3} on finite jet determination for all mappings between hypersurfaces of
D'Angelo finite type, we obtain the following .

\begin {Cor} Let $M$ and $M'$ be real-analytic hyperfaces in $\C^{n+1}$ of D'Angelo
finite type at $0$, with $(\2z,\2w)\in \bC^{n}\times\bC$ normal coordinates for
$M'$ at $0$.  Then there exists and integer $k>0$ with the following property.  Let
$H^1, H^2\colon  (\bC^{n+1},0)\to (\bC^{n+1},0)$ be germs at $0$ of holomorphic
mappings sending $M$ into
$M'$ and $H^j=(F^j,G^j)$, $j=1,2$, in the coordinates $(\2z,\2w)$. If $j^k_0
F^1=j^k_0F^2$, then $H^1\equiv H^2$.
\end{Cor}

\begin{proof}[Proof of Theorem $\ref{T:fd}$] Assume that the mappings in the class
$\mathcal F$ are determined by their $k_0$-jets at $0$. In view of Theorem
\ref{t:normal1}, it suffices to take normal coordinates $(\2z,\2w)\in
\bC^{n'}\times\bC^{d'}$ for $M'$ at $0$, write $H=(F,G)$ in these coordinates, and
show that there exists a $k\geq k_0$ such that the $k_0$-jet of $G$ at 
$0$ is determined by the $k$-jet of $F$ at $0$. We start with equation \eqref{e:Gdet},
which in view of Remark \eqref{r:form} can be written as follows
\begin{multline}\Label{e:Gdet2}
G(Z)=\Psi \big(\overline{F}\big(\hat h_1\big(\frac{Z}{\lambda^l},\lambda\big)\big)
,F\big(
\hat g_1\big(\frac{Z}{\lambda^l},\lambda\big)\big),\ldots, 
\\
\overline{F}\big(\hat
h_{d}\big(\frac{Z}{\lambda^l},\lambda\big)\big), F\big(\hat
g_d\big (\frac{Z}{\lambda^l},\lambda\big)\big),\overline{F}\big(\hat
h_{d+1}\big(\frac{Z}{\lambda^l},\lambda\big)\big),F(Z)\big),
\end{multline}
where the $\hat h_j$ and $\hat g_j$ are as in \eqref{e:hat}. By differentiating
\eqref{e:Gdet2} repeatedly with respect to $Z$, setting $Z=0$, we conclude that 
\begin{equation}\Label {e:series}
\frac{\partial^{|\alpha|} G}{\partial Z^\alpha}(0)=\sum_{j=-|\alpha|l}^\infty  a^\alpha_j
\lambda^j,
\end{equation}
where each coefficient $a^\alpha_j$ is a polynomial in the
components of $j^{|\alpha|(l+1)+j}_0F$ for all $j\geq -|\alpha|l$. Since the left hand 
side of \eqref{e:Gdet2} (and hence of \eqref{e:series})
is independent of $\lambda$, the coefficient $a^\alpha_j=0$ for $j\neq 0$ and $(\partial^
{|\alpha|} G/ \partial Z^\alpha) (0)=a^\alpha_0$. This completes the proof with
$k=k_0(l+1)$.
\end{proof}

\begin{Rem}{\rm The proof of Theorem \ref{T:fd}  shows that if the mappings $H$ in
$\mathcal F$ are determined by their $k_0$-jets at $0$, then they are also determined by
the $k$-jets at $0$ of $\2z\circ H$, where $k=k_0(l+1)$ and $l$ is the integer given
in Remark \ref{r:form} (depending only on $M$). However, we do not know any example
where the $k_0$-jet at $0$ of $\2z\circ H$ does not already determine $H$. If $M$
and $M'$ are strictly pseudoconvex hypersurfaces in $\bC^N$, then it follows from the
work of Chern--Moser \cite{cm} that $k_0$ can be taken to be $2$.  
Kruzhilin--Loboda
\cite{KL83} proved that for non-spherical strictly pseudoconvex hypersurfaces, the
stability group can be linearized in Chern--Moser normal coordinates and, hence, 
one may take $k_0=1$. In both the spherical and non-spherical case, one can check
directly that the $k_0$-jet of $\2z\circ H$ suffices to determine $H$. }
\end{Rem}

\end{document}